\documentclass[10pt]{article}

\usepackage{amscd,amsmath, amssymb, fancyhdr, mathbbol}

\numberwithin{equation}{section}

\def\eqref#1{(\ref{#1})}

\newcommand{\arrow}{{\:\longrightarrow\:}}
\newcommand{\Z}{{\Bbb Z}}
\def\C{{\Bbb C}}

\def\1{\sqrt{-1}\:}

\newcommand{\cntrct}                
{\hspace{2pt}\raisebox{1pt}{\text{$\lrcorner$}}\hspace{2pt}}

\newcommand{\calo}{{\cal O}}

\renewcommand{\bar}{\overline}
\renewcommand{\phi}{\varphi}
\renewcommand{\epsilon}{\varepsilon}
\renewcommand{\geq}{\geqslant}
\renewcommand{\leq}{\leqslant}

\newcommand{\Teich}{\operatorname{Teich}}

\newcommand{\Aut}{\operatorname{Aut}}

\newcommand{\Diff}{\operatorname{Diff}}

\newcommand{\Comp}{\operatorname{Comp}}

\newcounter{Mycounter}[section]
\newcounter{lemma}[section]
\renewcommand{\thelemma}{{Lemma \thesection.\arabic{lemma}}}
\newcommand{\lemma}{%
    \setcounter{lemma}{\value{Mycounter}}
    \refstepcounter{lemma}
    \stepcounter{Mycounter}
    {\noindent \bf \thelemma:\ }}

\newcounter{claim}[section]
\newcounter{sublemma}[section]
\newcounter{corollary}[section]
\renewcommand{\thecorollary}{{Corollary \thesection.\arabic{corollary}}}
\newcommand{\corollary}{%
    \setcounter{corollary}{\value{Mycounter}}
    \refstepcounter{corollary}
    \stepcounter{Mycounter}
    {\noindent \bf \thecorollary:\ }}

\newcounter{theorem}[section]
\renewcommand{\thetheorem}{{Theorem \thesection.\arabic{theorem}}}
\newcommand{\theorem}{%
    \setcounter{theorem}{\value{Mycounter}}
    \refstepcounter{theorem}
    \stepcounter{Mycounter}
    {\noindent \bf \thetheorem:\ }}

\newcounter{conjecture}[section]
\renewcommand{\theconjecture}{{Conjecture \thesection.\arabic{conjecture}}}
\newcommand{\conjecture}{%
    \setcounter{conjecture}{\value{Mycounter}}
    \refstepcounter{conjecture}
    \stepcounter{Mycounter}
    {\noindent \bf \theconjecture:\ }}

\newcounter{proposition}[section]
\newcounter{definition}[section]
\renewcommand{\thedefinition}
      {{Definition~\thesection.\arabic{definition}}}
\newcommand{\definition}{%
    \setcounter{definition}{\value{Mycounter}}
    \refstepcounter{definition}
    \stepcounter{Mycounter}
    {\noindent \bf \thedefinition:\ }}

\newcounter{example}[section]
\newcounter{remark}[section]
\renewcommand{\theremark}{{Remark \thesection.\arabic{remark}}}
\newcommand{\remark}{%
    \setcounter{remark}{\value{Mycounter}}
    \refstepcounter{remark}
    \stepcounter{Mycounter}
    {\noindent \bf \theremark:\ }}

\newcounter{problem}[section]
\newcounter{question}[section]
\makeatletter

\@addtoreset{equation}{section} \@addtoreset{footnote}{section}
\makeatother

\def\blacksquare{\hbox{\vrule width 5pt height 5pt depth 0pt}}
\def\endproof{\blacksquare}

\begin{document}
\begin{center}
{\LARGE\bf
Finiteness of Lagrangian fibrations with fixed invariants
\\[4mm]
}

Ljudmila Kamenova

\end{center}

{\small \hspace{0.1\linewidth}
\begin{minipage}[t]{0.8\linewidth}
{\bf Abstract} \\
We prove finiteness of hyperk\"ahler Lagrangian fibrations in any fixed 
dimension with fixed Fujiki constant and discriminant of the 
Beauville-Bogomolov-Fujiki lattice, up to deformation. 
We also prove finiteness of hyperk\"ahler Lagrangian fibrations with 
a very ample line bundle of a given degree on the general fiber of the 
fibration. 
\end{minipage}
}

{\scriptsize
\tableofcontents
}


\section{Introduction}


For a hyperk\"ahler manifold $M$, the Fujiki constant and the 
discriminant of the Beauville-Bogomolov-Fujiki lattice are topological 
invariants. It is very natural to fix them and ask for finiteness of 
hyperk\"ahler manifolds with these invariants. In this paper we establish 
finiteness of Lagrangian fibrations of hyperk\"ahler manifolds with fixed 
topological invariants as above. 

\hfill 

\theorem
There are at most finitely many deformation classes of Lagrangian fibrations 
$\pi: M \rightarrow \C P^n$ with a fixed Fujiki constant $c$ and a given 
discriminant of the Beauville-Bogomolov-Fujiki lattice $(\Lambda,q)$. 

\hfill

Francois Charles has the following boundedness result for families of 
hyperk\"ahler varieties up to deformation. He drops the assumption that 
$L$ is ample in Koll\'ar-Matsusaka's theorem applied for hyperk\"ahler 
manifolds and replaces it with the assumption that $q(L)>0$. 

\hfill 

\begin{theorem}(Charles, \cite{_Charles_})
Let $n$ and $r$ be two positive integers. Then there exists a scheme $S$ of 
finite type over $\C$, and a projective morphism ${\cal M} \arrow S$ such 
that if $M$ is a complex hyperk\"ahler variety of dimension $2n$ and $L$ is a 
line bundle on $M$ with $c_1(L)^{2n} = r$ and $q(L) > 0$, where $q$ is the
Beauville-Bogomolov form, then there exists a complex point $s$ of $S$ 
such that ${\cal M}_s$ is birational to $M$. 
\end{theorem}

\hfill

In our case, there is a natural line bundle $L$ associated to the 
Lagrangian fibration. Using Fujiki's formula, it is a straightforward 
observation that $q(L)=0$, while F. Charles deals with the case when 
$q(L)>0$ (in which case $M$ is projective by a result of D. Huybrechts: 
Theorem 3.11 in \cite{_Huybrechts:basic_}). 

\hfill

In the proof of our main theorems we use F. Charles' finiteness result 
applied to an ample line bundle with minimal positive square of the 
Beauville-Bogomolov-Fujiki form. Since we are interested in a finiteness 
result up to deformation equivalence, one can obtain an ample line bundle 
after deforming a given Lagrangian fibration to a projective one. We also 
use lattice theory estimates applied to the Beauville-Bogomolov-Fujiki form. 

\hfill

In \cite{_Sawon_finit_}, Sawon proved a finiteness theorem for Lagrangian 
fibrations with a lot of natural assumptions on the fibration, such as 
existence of a section, fixed polarization type of a very ample line bundle, 
semi-simple degenerations as the general singular fibers, and a maximal 
variation of the fibers. We give the precise statement of Sawon's theorem 
in \ref{C_S}. 
Due to a very recent progress of B. van Geemen and C. Voisin 
(\cite{_vG_V:matsushita_}) towards Matsushita's conjecture, the last 
condition of Sawon's theorem can be modified to only exclude 
isotrivial fibrations. 
Using the techniques in our proofs, one can also drop most of the 
other conditions in Sawon's theorem. We prove the following generalization. 

\hfill

\theorem 
Consider a Lagrangian fibration $\pi:M \arrow\C P^n$ such that 
there is a line bundle $P$ on $M$ with $q(P)>0$ and with a given $P$-degree 
$d$ on the general fiber $F$ of $\pi$, i.e., $P^n \cdot F = d$. 
Then there are at most finitely many deformation classes of hyperk\"ahler 
manifolds $M$ as above, i.e., they form a bounded family. 

\hfill 

For completeness of the exposition, we also mention Huybrechts' classical 
finiteness results. 

\hfill

\theorem (Huybrechts, \cite{_Huybrechts:finiteness_})
If the second integral cohomology $H^2(\Z)$ and the homogeneous polynomial
of degree $2n-2$ on $H^2(Z)$ defined by the first Pontrjagin class are 
given, then there exist at most finitely many diffeomorphism types of 
compact hyperk\"ahler manifolds of real dimension $4n$ realizing this structure.

\hfill

\theorem (Huybrechts, \cite{_Huybrechts:finiteness_}) \label{H2}
Let $M$ be a fixed compact manifold. Then there exist at most finitely
many different deformation types of irreducible holomorphic symplectic 
complex structures on $M$.

\hfill 

Using \ref{H2}, the author and Misha Verbitsky established the following 
finiteness results in \cite{_KV:fibrations_}. 

\hfill

\theorem (Kamenova-Verbitsky, \cite{_KV:fibrations_}) \label{KVfin}
Let $M$ be a fixed compact manifold. Then there are 
only finitely many deformation types of hyperk\"ahler 
Lagrangian fibrations $(M,I)\arrow \C P^n$, for all 
complex structures $I$ on $M$. 

\hfill

In the main theorem of this paper we prove finiteness of deformation 
classes of the total space $M$ of the Lagrangian fibration 
$M \arrow \C P^n$ with fixed dimension, Fujiki constant and discriminant 
of the Beauville-Bogomolov-Fujiki lattice. As a corollary of \ref{KVfin} 
one also obtains finiteness of the deformation classes of the Lagrangian 
fibration $M \arrow \C P^n$.


\section{Hyperk\"ahler geometry: preliminary results}


\subsection{Basic definitions}

\definition  A {\bf hyperk\"ahler manifold}
is a compact K\"ahler holomorphic symplectic manifold. 

\hfill

\definition
 A hyperk\"ahler manifold $M$ is called
{\bf simple} if $H^1(M, \C)=0$ and $H^{2,0}(M)=\C$.

\hfill

\theorem
(Bogomolov's Decomposition Theorem,
\cite{_Bogomolov:decompo_}, \cite{_Besse:Einst_Manifo_}). 
Any hyperk\"ahler manifold admits a finite covering,
which is a product of a torus and a finite collection of 
simple hyperk\"ahler manifolds. 
\endproof

\hfill

\remark
From now on, we assume that all hyperk\"ahler manifolds are simple. 

\hfill

\remark 
The following two notions are equivalent: 
a holomorphic symplectic K\"ahler manifold and 
a manifold with a {\em hyperk\"ahler structure}, 
that is, a triple of complex structures 
satisfying the quaternionic relations and 
parallel with respect to the Levi-Civita 
connection. In the compact case the equivalence 
between these two notions is provided by Yau's solution 
of Calabi's conjecture (\cite{_Besse:Einst_Manifo_}). 
In this paper we assume compactness and we use the 
complex algebraic point of view. 

\hfill

\definition
Let $M$ be a compact complex manifold and $\Diff^0(M)$ the connected 
component of the identity of its diffeomorphism group. 
Denote by $\Comp$ the space of complex structures on $M$, equipped with 
a structure of Fr\'echet manifold. The {\em Teichm\"uller space} of $M$ 
is the quotient  $\Teich:=\Comp/\Diff^0(M)$. 
For a hyperk\"ahler manifold $M$, the Teichm\"uller space is 
finite-dimensional (\cite{_Catanese:moduli_}). 
Let $\Diff^+(M)$ be the group of orientable diffeomorphisms of 
a complex manifold $M$. The {\em mapping class group} 
$$\Gamma:=\Diff^+(M)/\Diff^0(M)$$ acts naturally on $\Teich$. 
For $I\in \Teich$, let $\Gamma_I$ be the subgroup of $\Gamma$ which fixes 
the connected component of complex structure $I$. 
The {\em monodromy group} is the image of $\Gamma_I$ in $\Aut H^2(M, \Z)$. 

\subsection{The Beauville-Bogomolov-Fujiki form}

\theorem
(Fujiki, \cite{_Fujiki:HK_}) \label{Fujiki_formula}
Let $\eta\in H^2(M)$, and $\dim M=2n$, where $M$ is
hyperk\"ahler. Then $\int_M \eta^{2n}= c \cdot q(\eta,\eta)^n$,
for some integral quadratic form $q$ on $H^2(M)$, where $c>0$ is a 
constant depending on the topological type of $M$. The constant $c$ in Fujiki's 
formula is called the {\bf Fujiki constant}. 
\endproof

\hfill

\definition
This form is called the 
{\bf  Beauville-Bogomolov-Fujiki form}. 

\hfill

\remark
The form $q$ has signature $(3,b_2-3)$.
It is negative definite on primitive forms, and positive
definite on the space $\langle \Omega, \bar \Omega, \omega\rangle$
where $\Omega$ is the holomorphic symplectic form 
and $\omega$ is a K\"ahler form 
(see e. g. \cite{_Verbitsky:cohomo_}, Theorem 6.1,
or \cite{_Huybrechts:lec_}, Corollary 23.9).

\hfill

\definition
Let $[\eta]\in H^{1,1}(M)$ be a real (1,1)-class on
a hyperk\"ahler manifold $M$. We say that $[\eta]$
is {\bf parabolic} if $q([\eta],[\eta])=0$.
A line bundle $L$ is called {\bf parabolic} if the class $c_1(L)$
is parabolic.

\hfill

\remark \label{_P_L_identity_}
If $L$ is a parabolic class and $P \in H^2(M)$ is any class, then 
after we substitute $\eta = P+tL$ into Fujiki's formula in 
\ref{Fujiki_formula}, and compare the coefficients of $t^n$ on both sides, 
we obtain ${2n \choose n} P^n L^n = c 2^n q(P,L)^n$.

\subsection{The SYZ conjecture and Matsushita's conjecture}

\theorem
(Matsushita, \cite{_Matsushita:fibred_}).
Let $\pi:\; M \rightarrow B$ be a surjective holomorphic map
from a hyperk\"ahler manifold $M$ to a base $B$, with $0<\dim B < \dim M$.
Then $\dim B = 1/2 \dim M$, and the fibers of $\pi$ are 
holomorphic Lagrangian (this means that the symplectic
form vanishes on the fibers).

\hfill

\definition Such a map is called
{\bf a holomorphic Lagrangian fibration}.

\hfill

\remark The base of $\pi$ is conjectured to be
rational. J.-M. Hwang (\cite{_Hwang:base_}) 
proved that $B\cong \C P^n$, if it is smooth.
D. Matsushita (\cite{_Matsushita:CP^n_}) 
proved that it has the same rational cohomology
as $\C P^n$, if it is smooth.

\hfill

\definition A line bundle $L$ is called 
{\bf semiample} if  $L^N$ is generated 
by its holomorphic sections which have  
no common zeros. 

\hfill

\remark From semiampleness 
it trivially follows that $L$ is nef. Indeed,
let $\pi:\; M \rightarrow {\Bbb P}H^0(L^N)^*$ be the standard
map. Since the sections of $L$ have no common zeros, $\pi$ is 
holomorphic. Then $L\cong \pi^* \calo(1)$, and the
curvature of $L$ is the pullback of a 
K\"ahler form on $\C P^n$. However, 
a nef bundle is not necessarily semiample 
(see e.g. \cite[Example 1.7]{_Demailly_Peternell_Schneider:nef_}).

\hfill

\remark \label{_main_rem}
Let $\pi:\; M \rightarrow B$ 
be a holomorphic Lagrangian fibration, and $\omega_B$
a K\"ahler class on $B$. Then $\eta:=\pi^*\omega_B$ is 
semiample and parabolic. The converse is also
true, by Matsushita's theorem:
if $L$ is semiample and parabolic, $L$ induces a Lagrangian
fibration. 

\hfill

\conjecture\label{_SYZ_conj_Conjecture_}
({\bf Hyperk\"ahler SYZ conjecture})
Let $L$ be a parabolic nef line bundle
on a hyperk\"ahler manifold. Then
$L$ is semiample.

\hfill

\remark
The SYZ conjecture can be seen as
a hyperk\"ahler version of the ``abundance conjecture''
(see e.g. \cite{_Demailly_Peternell_Schneider:ps-eff_}, 
2.7.2). 

\hfill

\conjecture\label{_Matsushita_Conjecture_}
({\bf Matsushita's conjecture}) 
Every holomorphic Lagrangian fibration $\pi:M\rightarrow\C P^n$ is 
either locally isotrivial or the fibers vary maximally in the 
moduli space of Abelian varieties ${\cal A}_n$. 

\hfill

\remark
This conjecture was introduced to the author in private communictions 
with J. Sawon, D. Matsushita and J.-M. Hwang. 

\hfill

B. van Geemen and C. Voisin recently proved a weeker version of 
Matsushita's conjecture. 

\hfill

\begin{theorem} (B. van Geemen, C. Voisin, \cite{_vG_V:matsushita_}) 
\label{gv-m}
Let $X$ be a projective hyperk\"ahler manifold of dimension $2n$ admitting a 
Lagrangian fibration $f:X\rightarrow B$, where $B$ is smooth. 
Assume $b_{2, tr}(X) = b_2(X) - \rho (X) \geq 5$. Then a very general deformation 
$(X',f',B')$ of the triple $(X,f,B)$ satisfies Matsushita's conjecture. 
\end{theorem}

\subsection{Charles' and Sawon's finiteness theorems} \label{C_S}

Our main results in this paper rely on the following theorem. 

\hfill

\begin{theorem}(Charles, \cite{_Charles_}) \label{charles}
Let $n$ and $r$ be two positive integers. Then there exists a scheme $S$ of 
finite type over $\C$, and a projective morphism ${\cal M} \arrow S$ such 
that if $M$ is a complex hyperk\"ahler variety of dimension $2n$ and $L$ is a 
line bundle on $M$ with $c_1(L)^{2n} = r$ and $q(L) > 0$, where $q$ is the
Beauville-Bogomolov form, then there exists a complex point $s$ of $S$ 
such that ${\cal M}_s$ is birational to $M$. 
\end{theorem}

\hfill

We would also like to mention the following theorem in the recent literature. 

\hfill

\begin{theorem}(Sawon, \cite{_Sawon_finit_})
\label{_Sawon_finiteness_}
Fix positive integers $n$ and $d_1,\cdots,d_n$, with $d_1|d_2|\cdots |d_n$. 
Consider Lagrangian fibrations $\pi:M\rightarrow\C P^n$ that satisfy: 

(1) $\pi:M\rightarrow\C P^n$ admits a global section,

(2) there is a very ample line bundle on $M$ which gives a polarization of 
type $(d_1,\cdots,d_n)$ when restricted to a generic smooth fibre $M_t$,

(3) over a generic point $t$ of the discriminant locus the fibre 
$M_t$ is a rank-one semi-stable degeneration of abelian varieties,

(4) a neighbourhood $U$ of a generic point $t\in\C P^n$ describes a maximal 
variation of abelian varieties.

Then there are finitely many such Lagrangian fibrations up to deformation.
\end{theorem}

\hfill

\remark
Notice that as a corollary of Matsushita's conjecture, part $(4)$ of 
Sawon's Theorem simply excludes locally isotrivial fibrations. We need to 
apply only the deformational version (van Geemen-Voisin's \ref{gv-m}) 
of Matsushita's conjecture to Sawon's 
theorem in order to remove the seemingly restrictive assumption $(4)$. 

\hfill 

\remark
We would like to point out that 
if there is a section $\sigma : \C P^n \arrow M$, this means 
that $\sigma (\C P^n)$ would be a Lagrangian subvariety in $M$. 
Finding Lagrangian $\C P^n$'s in a hyperk\"ahler manifold is itself a 
very interesting task (for example, see \cite{_H_T_}). Moreover, 
the Lagrangian $\sigma (\C P^n)$ would have to intersect the general 
fiber of $\pi$ in one point.


\section{Main results}


Consider a lattice $\Lambda$, i.e., a free $\Z-$module of finite rank 
equipped with a non-degenerate symmetric bilinear from $q$ with 
values in $\Z$. If $\{ e_i \}$ is a basis of $\Lambda$, the 
{\it discriminant} of $\Lambda$ is defined as $\text{discr}(\Lambda) = 
\text{det}(e_i \cdot e_j)$. 

\hfill

\lemma\label{_bounded_lemma_}
Let $(\Lambda, q)$ be an indefinite lattice and $v \in \Lambda$ be an 
isotropic non-zero vector. Then there exists a positive vector 
$w \in \Lambda$ such that $0 < q(w,v) \leq |\text{discr}(\Lambda)|$ 
and $0<q(w,w) \leq 2 |\text{discr}(\Lambda)|$. 

\hfill

\begin{proof}
Let $w_0$ be a vector with minimal positive intersection $q(w_0,v)$. 
Then by Lemma 3.7. in \cite{_KV:fibrations_}, $q(w_0, v)$ divides 
$N = |\text{discr}(\Lambda)|$. Therefore, $0 < q(w_0, v) \leq N$. 
If $q(w_0,w_0) > 0$, let $\alpha$ be the smallest integer such that 
$q(w_0 + \alpha v, w_0 + \alpha v) >0$. Then we can take $w=w_0+\alpha v$. 
Otherwise, if $q(w_0,w_0) \leq 0$, 
consider the vectors $\{ w_0 + \alpha v \}$. Since $q(v,v)=0$, the square 
of such a vector is: $q(w_0 + \alpha v, w_0 + \alpha v) = q(w_0,w_0) + 
2 \alpha q(w_0, v)$. Take $\alpha$ to be a positive integer such that 
$q(w_0 + \alpha v, w_0 + \alpha v) > 0$. Set $w = w_0 + \alpha v$. 
Then in both cases $w$ is a positive vector with 
$0 < q(w,v) = q(w_0,v) \leq N$. 
Notice that automatically $0<q(w,w) = q(w_0 + \alpha v, w_0 + \alpha v) = 
q(w_0,w_0) + 2 \alpha q(w_0, v) \leq 2N = 2 |\text{discr}(\Lambda)|$. 
\end{proof}

\hfill






We recall the following result from a paper of the author's 
together with Misha Verbitsky, Theorem 3.6 in \cite{_KV:fibrations_}. 

\hfill

\theorem \label{KV_finit}
Consider the action of the monodromy group
$\Gamma_I$ on $H^2(M,\Z)$, and let $S\subset H^2(M,\Z)$
be the set of all classes which are parabolic
and primitive. Then there are only
finitely many orbits of $\Gamma_I$ on $S$.

\hfill

Our main result is the following finiteness theorem. 

\hfill

\theorem\label{_finiteness_Theorem_1_}
There are at most finitely many deformation classes of Lagrangian fibrations 
$\pi: M \rightarrow \C P^n$ with a fixed Fujiki constant $c$ and a given 
discriminant of the Beauville-Bogomolov-Fujiki lattice $(\Lambda,q)$. 

\hfill

\begin{proof}
As in \ref{_main_rem}, Lagrangian fibrations correspond to parabolic 
semiample classes. Now consider $S\subset H^2(M,\Z)$ defined 
above, the set of all classes which are parabolic and primitive, which 
is possibly larger than the set of parabolic semiample classes. 
By \ref{KV_finit}, there are only finitely many orbits of the 
monodromy group $\Gamma_I$ on $S$.

Let $L$ be a nef parabolic class ($q(L)=0$) coming from the Lagrangian 
fibration. Deform the Lagrangian fibration preserving the fibration 
structure, i.e., preserving the class of $L$ to a projective 
hyperk\"ahler Lagrangian fibration. Since we are interested in finiteness 
results up to deformation, we are going to work in the projective setting. 
By Huybrechts result (Theorem 3.11 in \cite{_Huybrechts:basic_}), 
there exists a line bundle with positive square. 
Apply \ref{_bounded_lemma_} for $(\Lambda,q) = (H^2(X,\Z),q)$ and $v=L$. 
There exists a positive vector $w$ with 
 $0 < q(w,v) \leq |\text{discr}(\Lambda)| =N$. We could choose $w$ to be 
a vector with the smallest positive square $q(w,w)>0$. From the lemma 
we see that $0<q(w,w) \leq 2 |\text{discr}(\Lambda)|$, which is bounded since 
we consider a fixed discriminant. 

Now we can apply F. Charles' \ref{charles} to the case when the first Chern 
class is $w$, in which case, by Fujiki's formula, 
$0< r = w^{2n} = c \cdot q(w,w)^n \leq c \cdot (2|\text{discr}(\Lambda)|)^n$ 
is bounded. For each $r$ in this interval we obtain only finitely many 
deformation classes of the total space $M$. 
\end{proof}

\hfill

Since the families of hyperk\"ahler manifolds as above 
form a bounded family, there are only finitely many 
choices of the second Betti number which plays an important role in 
studying the geometry of hyperk\"ahler manifolds. We obtain the following. 

\hfill

\corollary
In the assumptions of \ref{_finiteness_Theorem_1_}, the second Betti number 
$b_2(M)$ is bounded. 

\hfill

Using similar methods as above together with F. Charles' \ref{charles},  
we generalize Sawon's \ref{_Sawon_finiteness_} by dropping most of 
the assumptions. 

\hfill

\theorem \label{_finiteness_Theorem_2_}
Consider a Lagrangian fibration $\pi:M \arrow\C P^n$ such that 
there is a line bundle $P$ on $M$ with $q(P)>0$ and with a given $P$-degree 
$d$ on the general fiber $F$ of $\pi$, i.e., $P^n \cdot F = d$. 
Then there are at most finitely many deformation classes of hyperk\"ahler 
manifolds $M$ as above, i.e., they form a bounded family. 

\hfill

\begin{proof}
Let $L$ be a nef parabolic class ($q(L)=0$) coming from the Lagrangian 
fibration (e.g., as the pullback of a hyperplane class on $\C P^n$). 
The fundamental class $[F]$ of the general fiber of $\pi$ is proportional 
to $L^n$. We can fix the constant multiple in such a way that $[F] = L^n$. 
By assumption, $P^n \cdot L^n = d$ is fixed. 
Consider the classes $\{ P-kL \}$ for $k \in \Z$. We would like to bound 
the top degree of one of these classes, and apply F. Charles' theorem. 
Choose an integer $k \geq \frac{q(P,P) - 2 q(P,L)}{2 q(P,L)}$. For such $k$ 
we have the following estimate, where $c$ is the Fujiki constant: 
$(P-kL)^{2n} = c \cdot q(P-kL, P-kL)^n = c (q(P,P) - 2k q(P,L))^n \leq 
c (2 q(P,L))^n = 2^n c \cdot q(P,L)^n = {2n \choose n} P^n \cdot L^n = 
{2n \choose n} d.$ Here we applied Fujiki's formula twice 
(as in \ref{Fujiki_formula} and \ref{_P_L_identity_}). 
In order to apply \ref{charles}, we also need $q(P-kL, P-kL)>0$, i.e., 
$q(P,P)-2kq(P,L)>0$. Combining with the previous restriction on $k$, we have to 
choose $$k \in \Big[ \frac{q(P,P) - 2 q(P,L)}{2 q(P,L)}, 
\frac{q(P,P)}{2 q(P,L)} \Big).$$ 
Since $P$ is in the interior of the of the positive cone 
${\cal C}$ 
and $L$ is on the boundary of $\cal C$, it follows that $q(P, L)>0$ 
(Corollary 7.2 in \cite{_BHPV_}). 
The interval above is well-defined, because $q(P, L)>0$. 
For such a choice of the integer $k$, the top intersection of $P -kL$ is 
bounded and $q(P-kL) >0$. We apply F. Charles' \ref{charles} to obtain 
a bounded family of such $M$, which implies finiteness of deformations of $M$. 
\end{proof}

\hfill

\remark
In \ref{_finiteness_Theorem_1_} and \ref{_finiteness_Theorem_2_} we prove 
finiteness of deformation classes of the total space $M$ of the Lagrangian 
fibration. However, in \ref{KVfin} the author together with Misha Verbitsky 
prove that for a fixed compact manifold $M$ there are only finitely many 
deformation types of hyperk\"ahler Lagrangian fibrations with total space $M$. 

\hfill

\noindent{\bf Acknowledgments.} 
The author would like to express special thanks to Misha Verbitsky whose 
support on this project was invaluable. 
The author is grateful to Francois Charles for very interesting discussions 
on this subject and for his help on the exposition. 

\hfill

{\small

\hfill

\noindent {\sc Ljudmila Kamenova\\
Department of Mathematics, 3-115 \\
Stony Brook University \\
Stony Brook, NY 11794-3651, USA \\
 }

\end{document}